\documentclass[11pt]{article}
\usepackage{amsmath}
\usepackage{amsfonts}
\usepackage{amssymb}
\usepackage{amsthm}

\def\ball{{\cal B}}

\def\sphere{{\cal{S}}}
\def\ball{{\cal B}}

\def\RR{{\bf R}}
\def\NN{{\bf N}}
\def\ZZ{{\bf Z}}
\def\con{{\bf C}}

\def\proof{{\bf Proof: \, }}
\def\beqns{\begin{eqnarray*}}
\def\eeqns{\end{eqnarray*}}
\def\beqn{\begin{eqnarray}}
\def\eeqn{\end{eqnarray}}

\def\laplace{\Delta}
\def\dlaplace{\overline{\Delta}}
\def\domain{{\Omega}}

\def\supp{{\rm supp} \,}

\def\eop{{\vrule height7pt width7pt depth0pt}\par\bigskip}
\def\meop{\qquad \mbox{\eop}}

\newtheorem{theorem}{Theorem}
\newtheorem{proposition}[theorem]{Proposition}
\newtheorem{definition}[theorem]{Definition}
\newtheorem{lemma}[theorem]{Lemma}
\newtheorem{corollary}[theorem]{Corollary}
\newtheorem{assumption}[theorem]{Assumption}
\newtheorem{example}{Example}
\theoremstyle{remark}

\newcommand{\spam}{\mathop{\mathrm{span}}}

\newcommand{\dif}{\mathrm{d}}

\newcommand{\rbf}{\widetilde{\phi}_{d,k}}

\renewcommand{\d}{\mathrm{dist}}

\newcommand{\B}{\mathcal{B}}

\begin{document}

\title{On the density of polyharmonic splines  
      }

\author{ T.~Hangelbroek\thanks{Corresponding author. E-mail: {\tt hangelbr@math.hawaii.edu} }
\thanks{Department of Mathematics, University of Hawaii, 
 Honolulu, HI 96822, USA.   Research supported by  by grant DMS-1232409 from the National
    Science Foundation.},
    J.~Levesley\thanks{Department of Mathematics, University of Leicester, Leicester LE1 7RH , UK.}} 

\maketitle
\begin{abstract}
This article treats the question of fundamentality of the translates of a polyharmonic spline kernel (also known as 
a surface spline) in the space of continuous functions on a compact set $\Omega\subset \RR^d$
when the translates are restricted to $\Omega$. Fundamentality is not hard to demonstrate
 when a low degree polynomial may be added
or when translates are permitted to lie outside of $\Omega$;
the challenge of this problem stems from the presence of the boundary,
for which all successful approximation schemes require an added polynomial. 

When $\Omega$ is the unit ball, we demonstrate that translates of polyharmonic splines are fundamental
by considering two related problems: the fundamentality in the space of  functions vanishing at the boundary
and fundamentality of the restricted kernel in the space of continuous function on the sphere.
This gives rise to a new approximation scheme composed of two parts: one which approximates purely on $\partial \Omega$, and a second part involving a shift invariant approximant of a function vanishing outside
of a neighborhood $\Omega$.
\end{abstract}
\section{Introduction}

Interpolation using the thin-plate spline radial basic function $\phi(r)=r^2 \log r$ requires the addition of a linear polynomial.
This requirement is necessary to guarantee the solvability of the interpolation problem (see e.g. Light and Wayne \cite{light2}).

It can also be shown that for any compact set $\domain \subset \RR^d$ 
functions of the form
$$
\sum_{y \in Y} \alpha_y \phi(|x-y|)+p_1(x),
$$
where $Y \subset \domain$ is a finite set and $|\cdot  |$ is the Euclidean norm, 
are dense in $\con(\domain)$ -- the space of continuous functions. Here $p_1$ is a polynomial of degree $1$. 
For general compact sets $\Omega$, density requires the addition of an extra (polynomial) term. 
This is clearly the case when $\Omega = \{x\}$ is a singleton. 
A more interesting example is $\Omega  = \sphere^1 \cup \{0\}\subset \RR^2$,  the unit circle augmented with the origin. 
In this case, if $Y\subset \Omega$,
then any combination 
$\sum_{y \in Y} \alpha_y \phi(|\cdot-y|)$
 vanishes at the origin.

In this article, we are concerned with density of translates of {\em radial basis functions} $\phi$,
a popular tool used to construct approximants on $\RR^d$  . 
This problem goes by the name of  {\em fundamentality}.
\begin{definition}
A set $A$ is fundamental in a normed vector space $X$ if, for every $x
\in X$, given $\epsilon > 0$, there exists a finite subset $Y
\subset A$ and coefficients $a_y$, $y \in Y$, such that
$$
\left\| x - \sum_{y \in Y} a_y y \right\| < \epsilon.
$$
\end{definition}
Given a function $\phi:[0,\infty)\to \RR$, the basic problem is to find 
compact sets $\Omega\subset \RR^d$
for which 
$A(\Omega,\phi):=\{ \phi (|\cdot - y|)\mid y\in \Omega\}$
is fundamental 
 in the space of continuous functions $\con(\Omega)$ equipped with the uniform norm
$\|  \cdot \|_{\infty,\domain}$.
We consider this problem for the family of polyharmonic splines %
\footnote{So called, because $\phi_{d,k}$ is, up to a constant factor, the fundamental solution for the $k$ fold Laplacian in $\RR^d$} 
\begin{equation}\label{PHS_definition}
\phi_{d,k}(r)= \left \{ \begin{array}{cc}
r^{2k-d}, & d \;\; {\rm odd}, \\
r^{2k-d} \log r, & d \;\; {\rm even},
\end{array} \right.
\end{equation}
 where $k>d/2$ and $d$ is the spatial dimension. 
 As with the case of  thin plate splines ($k=2$, $d=2$), density is guaranteed with the addition of certain low order polynomial terms. However, little is known when approximating only by translates of the polyharmonic splines. 
 
 In Cheney and Light \cite{cheney} there is the excellent Chapter 18 which discusses the fundamentality of translates of positive definite and conditionally positive definite functions of order 1. This is a general setup that provides density for many radial basis functions, but not the polyharmonic splines in dimension $d\ge2$ (these are conditionally positive definite of order $k>d/2\ge 1$).
The arguments presented therein follow that of Brown \cite{brown}, which are measure-theoretic and non-constructive. 
There is a measure theoretic proof for density of the thin-plate spline $\phi(r)=r^2\log r$ in $\RR^2$, shown to us by Johnson \cite{johnson}, but this approach does not generalize to other space dimensions.

Fundamentality has been studied by Jackson 
in a setting very close to the one we consider.
He considers radial basis functions $\phi$ that are homogeneous (this includes $\phi_{d,k}$ in the case that $d$ is odd),
but he enlarges the fundamental set $A$ in a critical way.
Namely,  \cite[Theorem 4]{Jackson}  establishes the
fundamentality of  $A(\widetilde{\Omega},\phi)$
in $\con(\Omega)$
for compact $\Omega$ 
where $\widetilde{\Omega}\supset  \Omega$ is open.
Notably, no extra polynomial terms are required.

At this point, a key aspect of the problem emerges. 
When the translates $Y$ are forced to lie within the set $\Omega$, 
the question of fundamentality is a more challenging one. 
Indeed, similar results to those of Jackson 
can be obtained for all polyharmonic kernels 
$\phi_{d,k}$   with $k>d/2$ 
(for $d$ of any parity).
That is, for any open neighborhood $\widetilde{\Omega}\supset \Omega$,
 the set $A(\widetilde{\Omega}, \phi_{d,k})$  is fundamental in $\con(\Omega)$. 
This is  a consequence of Theorem \ref{main_general} below,
although it follows in a similar way from several earlier works. 

The consideration of $A(\widetilde{\Omega}, \phi_{d,k})$ is a key relaxation which avoids {\em boundary effects},
 where translates $Y$ are kept to one side of a boundary $\partial \Omega$.
The effect of a boundary  in radial basis function approximation is known to be detrimental in a number of ways:
it curtails approximation rates, forces instability by way of increased Lebesgue constants, causes
blow up of error near the boundary, etc. 
The fundamentality of the thin plate spline $\phi_{2,2}$
in the closed unit disk in $\RR^2$,  denoted  $\ball^2$, 
has been demonstrated in \cite{Hang}.
Specifically, \cite[Theorem 10]{Hang} 
establishes that 
 $A(\ball^2,\phi_{2,2})$ is fundamental in $\con(\ball^2)$
 by developing an approximation scheme $T_{Y} f = \sum_{y \in Y} \alpha_y \phi_{2,2}(|\cdot-y|)$
 that converges on  $\con^4(\ball^2)$ as $Y\subset \ball^2$ becomes dense. 
Although this scheme has been generalized in \cite{HangThesis}
 to other compact planar domains $\Omega\subset\RR^2$, these generalizations require addition of a linear polynomial.
 To date analogous results for other spatial dimensions ($d>2$) and other polyharmonic splines $\phi_{d,k}$, $k>d/2$,
 are elusive.
We pose this question because we wish to know when (or if) a polynomial term is needed to construct a good approximation and, secondly, in the hope that it will
generate approximation schemes
that adequately treat the boundary effects.

 In this paper we will show that translates of the polyharmonic splines $\phi_{d,k}$ with $k>d/2$
are fundamental in the space of continuous functions on the closed unit ball $\ball^d = \{ x \in \RR^d : | x | \le 1 \}$. 
As a step on the way we shall prove that for any compact set $\Omega\subset \RR^d$, such functions are fundamental in the set 
$\con_0(\Omega)$, of continuous functions which vanish on the boundary. Consequently, the
solution of the problem of density in $\con(\Omega)$ hinges upon  fundamentality on the boundary, in other words,
density of $\spam_{y\in\partial \Omega} \left.\bigl(\phi_{d,k}(\cdot -y)\bigr)\right|_{\partial\Omega}$ in $\con(\partial \Omega)$.

As a final note, we mention that for even spatial dimensions $d\in 2\NN$, the polyharmonic spline $\phi_{d,k}$ 
is not homogeneous. 
Consequently, the underlying spaces  are not dilation invariant: 
$$\text{for $\rho>0$,}\quad\spam_{y\in\rho\Omega} \phi_{d,k}\left(\cdot -y\right) \ne \spam_{y\in \Omega} \phi_{d,k} \left(\frac{\cdot}{\rho} - y\right).$$
Consequently, fundamentality of $A(\Omega,\phi_{d,k})$ in $\con(\Omega)$ does not necessarily imply fundamentality of
 $A(\rho\Omega,\phi_{d,k})$ in  $\con(\rho \Omega)$.
We investigate this in detail in Section \ref{arbitrary_radius} -- in particular we show that there are some spheres  
$\sphere^{d-1}(\rho) = \{x\in \RR^d\mid |x| = \rho\}$ 
for which density fails;  there exist values of $\rho>0$ so that 
$A(\sphere^{d-1}(\rho),\phi_{d,k})$ is not fundamental in $\con(\rho\sphere^{d-1})$.

\section{Density in $\con_0(\Omega)$, $\con(\partial \Omega)$ and $\con(\Omega)$}
In this section we discuss approximation by integer shifts of a function satisfying a simplified type of Strang-Fix condition.
We use this scheme to approximate functions vanishing on the boundary of a compact set $\Omega$. We then
discuss the polyharmonic B-splines --  particularly that they satisfy such a condition. Finally, we demonstrate that the polyharmonic splines are fundamental in $\con(\Omega)$ if 
their restriction  to $\partial \Omega$ is
fundamental in $\con(\partial \Omega)$.

\subsection{Approximation in Euclidean space}
A standard way of constructing an approximation to a function is via convolution with an approximation to the identity. Given $f \in \con(\RR^d)$ and $g \in L_1(\RR^d)$, 
the convolution is
$$
f \ast g(x) = \int_{\RR^d} f(x-y) g(y) dy.
$$
If we further require that $\int_{\RR^d} g =1$, then we can construct the approximate identity $g_h(x) = h^{-d} g(h^{-1} x)$, $x \in \RR^d$, $h > 0$. 
In this case, for $f\in\con_0(\Omega)$, we have that $g_h*f \to f$ uniformly  as $h\to 0$.

If we discretize the convolution $f \ast g_h$ on the $h$-scaled multi-integer grid (this is sometimes
referred to as the {\em semi-discrete convolution}) we arrive at the quasi-interpolant
\beqns
s_h (f,x) & = & \sum_{z \in \ZZ^d} f(hz) g(h^{-1}x - z). \label{quasi}
\eeqns

We apply $s_h$ to  functions that are uniformly continuous on $\RR^d$, i.e. functions for which the modulus of continuity 
$\omega(f,\delta,\Omega)  :=  \sup \{ f(x)-f(x-y) : x \in \Omega, \; \| y\|
\le \delta \}$
(suppressing $\Omega$ when $\Omega=\RR^d$) satisfies
$$
\omega(f,\delta) = \omega(f,\delta,\RR^d)  \rightarrow 0,
$$
as $\delta \rightarrow 0$. 
We will show that the quasi-interpolant converges uniformly on uniformly continuous functions
provided the following three conditions are satisfied.
\begin{assumption}\label{assumption}We assume the function $g\in \con(\RR^d)$ satisfies the following: 
\begin{enumerate}
\item[(I)] There exists an $M$ such that for every $x \in \RR^d$,
$
\sum_{z \in \ZZ^d} |g(x-z)| < M.
$

\item[(II)] The translates of $g$ form a partition of unity: i.e., for $x \in \RR^d$,
$
\sum_{z \in \ZZ^d} g(x-z) = 1.
$
\item[(III)] There exists an $M'$ such that for every $x \in \RR^d$,
$
\sum_{z \in \ZZ^d} |x-z||g(x-z)| < M'.
$
\end{enumerate}
\end{assumption}
We note that Assumption \ref{assumption} {\em(I)} implies that the $L_{\infty}$ operator norm of the semi-discrete convolution is bounded by $M$. This is well known (cf. \cite[Theorem 2.1]{JiaMic}) and easy to prove. (A variant of this is demonstrated below in Lemma \ref{little}.)

The rate of convergence of the quasi-interpolant operator $s_h$  has been studied extensively under the umbrella of {\em shift invariant approximation}. 
The approximation power of  $s_h$,
as well as other shift-invariant approximation schemes, 
is related to a Strang -- Fix condition on the basic function $g$ 
(this is a condition on the behavior of the Fourier transform of $g$ at $2\pi\ZZ^d$).  
We note that Assumption \ref{assumption} guarantees
that $g$ satisfies a Strang -- Fix condition of order zero (i.e., $\widehat{g}$ vanishes at the $2\pi\ZZ^d\setminus \{0\}$ 
with a simple zero) because the semi-discrete convolution reproduces constant functions.
For this article, we do not require the full force  of this theory, but direct the interested reader
to  \cite{dBR} for an in depth discussion of approximation rates for this operator.
Instead, we develop what is necessary for our purposes in the following proposition. 
\begin{proposition} \label{euclid}
Let $f$ be bounded and uniformly continuous on $\RR^d$. 
Then, given $\epsilon > 0$, there exists an $H>0$ such that $\| f - s_h(f,\cdot) \|_\infty < \epsilon$ for all $h < H$.
\end{proposition}
\proof Since the translates of $g$ form a partition of unity (item {\em (II)} of Assumption \ref{assumption}), we have the bound
\beqns
|f(x)-s_h(f,x)| & = & \left | f(x) -\sum_{z \in \ZZ^d} f(hz) g(h^{-1} x - z) \right | \\
& \le & \sum_{z \in \ZZ^d} | f(x) - f(hz) | | g(h^{-1} x - z) |.
\eeqns
This series can be split into a proximate and remote part, using $B_\delta(x) = \{ y \in \RR^d : | x - y | \le \delta \}$,   
the ball centered at $x$ having radius $\delta>0$ (to be determined). We have
\beqns
|f(x)-s_h(f,x)|& \le & \sum_{hz \in B_\delta(x)} | f(x) - f(hz) | | g(h^{-1} x - z) |\\
&\mbox{}& +
 \sum_{hz \not \in B_\delta(x)} | f(x) - f(hz) | | g(h^{-1} x - z) | \\
& \le & \omega(f,\delta) \sum_{hz \in B_\delta(x)}  | g(h^{-1} x - z) | + 2 \|f \|_\infty \sum_{hz \not \in B_\delta(x)}  | g(h^{-1} x - z) |, \\
\eeqns 
We may use Assumption \ref{assumption} {\em(I)} to obtain
$\omega(f,\delta) \sum_{hz \in B_\delta(x)}  | g(h^{-1} x - z) |
\le 
M \omega(f,\delta)$.
In order to control the second term, we use H{\" o}lder's inequality and Assumption \ref{assumption} {\em(III)}  , obtaining
 the inequality
$$
\sum_{hz \not \in B_\delta(x)}  | g(h^{-1} x - z) |
 \le 
 \sup_{hz \not \in B_\delta(x)} |h^{-1} x - z|^{-1} \sum_{hz \not \in B_\delta(x)} |h^{-1} x - z| | g(h^{-1} x - z) | 
 \le \frac{h}{\delta} 
M'  .
 $$
From this, we have
$
 2 \|f \|_\infty \sum_{hz \not \in B_\delta(x)}  | g(h^{-1} x - z) |
  \le 
  2 \|f \|_\infty {h \over \delta} M' .
$

Finally, let us choose $\delta$ such that
$\omega(f,\delta) < \epsilon/(2M)$,
and following this choice of $\delta$, choose
$
h < \epsilon \delta  /(4 M' \| f \|_\infty).
$
It follows that
\beqns
|f(x)-s_h(f,x)| 
 \le  M \omega(f,\delta)+   2 \|f \|_\infty {h \over \delta} M'<\epsilon
\eeqns
holds.\meop

\subsection{Approximating functions vanishing on the boundary of a compact set}

The difficulties one faces when proving density results via discretisation of convolutions happen due to the boundaries of the region. Because of this it is much easier to prove results where the function is zero on the boundary of the domain of interest. This shall be the case for this section. Since our target function is continuous on the compact domain $\Omega$, it is uniformly continuous. If we let $\bar f$ to be defined on the whole of $\RR^d$ by setting its value to zero on the complement of $\Omega$, it is easy to show that $\bar f$ is also uniformly continuous on $\RR^d$.

We define our approximation
\begin{equation}\label{Approximant_Definition}
s_h^\Omega(f,x) := \sum_{hz \in \Omega\cap h\ZZ} f(hz) g(h^{-1} x-z), \quad x \in \RR^d.
\end{equation}

Our first result, which we will use later, says that if $f$ is small on $\Omega$ then $s_h^\Omega(f,\cdot)$ is small everywhere.
\begin{lemma} \label{little}
Suppose $f \in \con(\Omega)$ and $\| f \|_{\infty,\Omega} \le \epsilon$. Suppose further that $g$ satisfies item
(I) of Assumption \ref{assumption}. Then,  $\| s_h^\Omega(f,\cdot) \|_{\infty} \le M \epsilon$.
\end{lemma}
\proof Since $g$ satisfies Assumption \ref{assumption} item {\em (I)} we have, for $x \in \RR^d$,
\beqns
| s_h^\Omega(f,x) | & \le & \left | \sum_{hz \in \Omega} f(hz) g(h^{-1} x-z)
\right | 
\le \|f\|_{\infty,\Omega} \sum_{hz \in \Omega} | g(h^{-1} x-z) | 
\le M \epsilon. \meop
\eeqns

Since, for $x \in \Omega$, $s_h^\Omega(f,x) = s_h(\bar f,x)$ as a corollary of Proposition~\ref{euclid}, we get our first approximation result for a compact domain:
\begin{corollary} \label{zero}
Let $f \in \con_0(\Omega)$ and $g$ satisfy Assumption \ref{assumption}. Then, given $\epsilon > 0$ there exists an $H>0$ such that $\| f - s_h^\Omega(f,\cdot) \|_{\infty,\Omega} < \epsilon$ for all $h < H$.
\end{corollary}

\subsection{Polyharmonic B-splines}
We employ the fact that $\phi_{d,k}$ is (up to a constant factor) the fundamental solution of the $k$-fold Laplacian in $\RR^d$, i.e., in the distributional sense, 
$$
\laplace^k \phi_{d,k} = C_{d,k} \delta 
\quad
\text{ where } 
\quad
C_{d,k} := \Gamma(d/2) \left[2^k \pi^{d/2} (k-1)! \prod_{\substack{i=0\\ i\ne m-d/2}}^{m-1} (2k-2i-d)\right]^{-1}.
$$
 Let $\{e_j : j=1,\cdots,d\}$ be an orthonormal basis for $\RR^d$. We can now form the discrete Laplacian of a continuous function $f$
$$
\dlaplace_h f (x) = \frac{1}{h^{2}}\left[2d\, f(x) - \sum_{j=1}^d \bigl(f(x+he_j)+f(x-he_j)\bigr)\right].
$$
Then, as suggested by Rabut \cite{CR1}, we can form the {\sl polyharmonic B-splines}\footnote{Rabut's definition differs from this one by a factor of $h^d$.}
$$
B_{d,k,h}  := {1 \over C_{d,k}} \dlaplace_h^k \phi(|\cdot|).
$$
We note that this kind of localization of the polyharmonic splines has an important history in radial basis function approximation, having occurred earlier as an analytic tool in \cite{BSW} and generating preconditioners in \cite{DynLevinRippa}.

When $h=1$, we drop the third subscript: $B_{d,k} := B_{d,k,1}$
and $ \dlaplace^k := \dlaplace_1^k$.
It follows that $B_{d,k} =  C_{d,k}^{-1} \dlaplace^k \phi(|\cdot|)= \sum_{z \in \ZZ^d} c_{z} \phi(|x-z|)$ for a 
sequence
$\mathbf{c} =(c_{z})_{z \in \ZZ^d }$ supported in the set $\{z\in \ZZ^d\mid \|z\|_{\ell_1} \le k\}$ (in other words, $c_{z}= 0$ for all 
multi-integers $z$ satisfying 
$\sum_{j=1}^d |z_j|>k$).
 We obtain
$$B_{d,k,h} = {h^{-2k}}\sum_{z \in \ZZ^d} c_{z} \phi(|x-hz|)$$
by rescaling the coefficients.
\begin{proposition} [From \cite{CR1}]\label{rabutres}
The function $B_{d,k} = B_{d,k,1}$ satisfies Assumption \ref{assumption}.
Additionally, the functions satisfy the dilation condition: for any $h>0$,
$$
  B_{d,k,h}(x) = h^{-d}  B_{d,k}\left(\frac{x}{h}\right).
$$
\end{proposition}
\proof
Assumption \ref{assumption} follows from \cite{CR1}.  
Indeed, items  {\em (I)} and {\em (III)} are a consequence of the decay $|B_{d,k}(x)|\le C(1+|x|)^{-d-2}$, which is  \cite[Theorem 2(i)]{CR1}.
Item {\em (II)} is precisely \cite[Theorem 4(i)]{CR1}.

When $d$ is
even, $\phi(|hy|) = h^{2k-d} \phi(|y|) + (h^{2k-d} \log h ) \,|y|^{2k-d}$. We have also that the functional
$\sum c_{z} \delta_z$ annihilates $\Pi_{2k-1}$, so
 $$\sum_{z \in \ZZ^d}c_{z} \phi\left(\left| {hy-hz} \right|\right)   = 
 h^{2k-d}\sum_{z \in \ZZ^d}c_{z} \left[\phi\left(\left|{y-z}\right|\right)
 + 
\log h   |y-z|^{2k-d} \right]
=
h^{2k-d}\sum_{z \in \ZZ^d}c_{z} \phi\left(\left|{y-z}\right|\right). $$
Letting $y = x/h$, it follows that 
$$ B_{d,k,h}(x) = B_{d,k,h}(hy)  = {h^{-2k}} \sum_{z \in h\ZZ^d}c_{z} \phi(|hy-hz|)  = h^{-d} \sum_{z \in \ZZ^d}c_{z} \phi\left(\left|{y-z}\right|\right)
= h^{-d}B_{d,k} \left(\frac{x}{h}\right).$$
When $d$ is odd, $\phi_{d,k}$ is homogeneous, and the dilation condition follows by a similar
(in this case simpler) argument.  
\meop

Although we are primarily interested in proving density properties of $B_{d,k}$, there exists an extensive
literature on convergence rates for the operator $s_h$ and related schemes.  We note in particular that
the function $B_{d,k}$ satisfies a Strang--Fix condition of order $2m$, which can be observed directly from the 
Fourier transform of $B_{d,k}$ in \cite[Theorem 2]{CR1}.
We refer the reader to \cite[Theorem 3.6]{dBR} for  a discussion of the approximation orders of 
this kind of approximation scheme.

\subsection{Approximating on compact sets}
We now give our main result.
\begin{theorem}\label{main_general}
Let $\Omega\subset \RR^d$ be a compact set. Assume that 
$A(\partial\Omega,\phi_{d,k})$
is fundamental in $\con(\partial \Omega)$. Let
$f \in \con(\Omega)$. 
Given $\epsilon > 0$ 
there is a finite subset 
$Y \subset \Omega$ 
and a set of coefficients 
$\{ \alpha_y : y \in Y \}$ 
such that
$$\left|f(x) - \sum_{y \in Y} \alpha_y \phi_{d,k}(|x-y|)\right| \le \epsilon$$ 
for all $x \in \Omega$.
\end{theorem}
\proof
By assumption,
we can select a finite subset $W \subset \partial \Omega$ and coefficients $\{ \alpha_y : y \in W \}$ such that
$$
\left|f(x) - \sum_{y \in W} \alpha_y \phi_{d,k}(|x-y|)\right| 
\le 
{\epsilon \over 4M},
$$
for all $x \in \partial \Omega$. 
Here $M$ is the specific constant used in Assumption \ref{assumption} {\em (II)}.

For $\gamma>0$ let us define the following neighborhoods:
\begin{itemize}
\item a $\gamma$ neighborhood of $\Omega$: 
$$\Omega_\gamma 
:= 
\{ x \mid  \d(x,\Omega) \le \gamma\},
$$
\item an annular neighborhood of the boundary:
$$
A_\gamma 
:= 
\{ x \mid \d(x,\partial \Omega) \le \gamma \},
$$ 
\item a one-sided neighborhood
$$
S_\gamma 
:=  
\{ x \in \Omega \mid  \d(x,\partial\Omega)  \le \gamma\}.
$$
\end{itemize}
We can extend $f$ beyond $\Omega$. For $\gamma=1$, we take $f^{e}$ to be a continuous extension of $f$ to  $\Omega_1$. (The extension beyond $\Omega$ is important for the approximant we will construct, although the values from $\Omega_1\setminus \Omega$ are not explicitly used.)

Now let us set $r(x) = f^e(x) - \sum_{y \in W} \alpha_y \phi_{d,k}\bigl(|x-y|\bigr)$, $x \in \Omega_1$. Then, since $r$ is continuous and small on $\partial \Omega$, there is a $0<\delta<1$ such that
$$
|r(x)| \le {\epsilon \over 2M}
$$
for all $x \in A_\delta$. By multiplying by a continuous cut-off function $\sigma$
$$
\sigma(x) = 
\begin{cases} 1, \quad& x\in \Omega\\ 1- \frac{\d(x,\partial\Omega)}{\delta}, &0,<\d(x,\Omega)<\delta,\\
0,&\d(x,\Omega)\ge\delta,
\end{cases}
$$
there
 is a continuous function $\bar r  = \sigma \times r$ such that $\bar r = r$ on 
$\Omega$, 
$\bar r (x) = 0$ when $x \not \in \Omega_{\delta}$, and 
$\| r \|_{\infty,S_\delta}\le \|\bar r \|_{\infty,A_\delta} \le \epsilon/(2M)$. 

We are now in a position to construct an approximant to $r$ of the form (\ref{Approximant_Definition}).
To this end, we assume that the kernel $g$ of Assumption \ref{assumption} is $B_{d,k}$.
Thus, using Lemma~\ref{little} we have, for any $h>0$,
\beqn
\left\| s_h^{A_\delta}(\bar r) \right\|_\infty & \le & {\epsilon \over 2}. \label{adel}
\eeqn
Since $B_{d,k}$ satisfies Assumption \ref{assumption} 
there exists an $H>0$ such that, for $h<H$
\beqn
\left\| \bar r-s_h^{\Omega_{\delta}}(\bar r) \right\|_\infty & \le & {\epsilon \over 2}. \label{odel}
\eeqn
For $x\in \Omega$, we approximate $f(x)$  with 
$s(x)
:=
\sum_{y \in W} \alpha_y \phi_{d,k}\bigl(|x-y|\bigr)
+
s_h^{\Omega\setminus S_{\delta}}(r,x)
$. 
Then, since 
$s_h^{\Omega\setminus S_{\delta}}(r,x)=s_h^{\Omega_{\delta}}(\bar r,x)-s_h^{A_\delta}(\bar r,x)$,
\beqns
\| f - s\|_{\infty,\Omega} & = & \| f - \sum_{y \in W} \alpha_y \phi_{d,k}(|x-y|)- s_h^{\Omega\setminus S_{\delta}}(r) \|_{\infty,\Omega} \\
& \le & \| \bar r - s_h^{\Omega_{\delta}}(\bar r)+s_h^{A_\delta}(\bar r) \|_{\infty} 
\le \epsilon,
\eeqns
for $h<H$, using (\ref{adel}) and (\ref{odel}). 

We make the further requirement that $h<\delta/k$ (where $\delta$ was chosen above, and $k$ is the order of polyharmonicity of $\phi_{k,d}$).

Finally, by its definition, (\ref{Approximant_Definition}), the translates of  $g_h = B_{d,k}\bigl(\cdot/h \bigr)$ used to construct $s_h^{\Omega\setminus S_{\delta}}$ 
lie in $\Upsilon:= h\ZZ^d\cap (\Omega\setminus S_{\delta})$.
Indeed,
$$s_h^{\Omega\setminus S_{\delta}} (r)
= 
\sum_{\xi \in h\ZZ^d\cap \Omega\setminus S_{\delta}} 
r(\xi) B_{d,k}\bigl(h^{-1}(\cdot-\xi)\bigr).$$
Fix $\xi\in h\ZZ^d\cap (\Omega\setminus S_{\delta})$. The basic function $B_{d,k} \bigl(h^{-1}(\cdot-\xi)\bigr)$ is a linear combination of the $h\ZZ^d$ translates of radial basis function $\phi_{d,k}$. Indeed, by  Proposition~\ref{rabutres} we have 
$$
B_{d,k} \bigl(h^{-1}(\cdot-\xi)\bigr)
 = 
 h^d B_{d,k,h} \bigl(\cdot-\xi \bigr)
=
h^{d-2k} \sum_{z\in\ZZ^d} c_{z}  \phi_{d,k}\bigl(|\cdot-\xi - hz\bigr|).
 $$
 It follows that $B_{d,k} \bigl(h^{-1}(\cdot-\xi)\bigr)$ is a linear combination of the  shifts of $\phi_{d,k}$
 taken from the set $\xi +h \:\supp{\mathbf{c}}$ (where $\mathbf{c}$ is the sequence  $\mathbf{c} = (c_z)_{z\in\ZZ^d}$).

Because the support of the sequence $\mathbf{c}$ lies in the ball of radius $k$ centered at the origin,
 $\supp{\mathbf{c}} \subset \left\{z\in \ZZ^d \mid \sum_{j=1}^d |z_j| \le k\right\} \subset \ball^d(k)$, and because $kh < \delta$,
 it follows that  the function
 $s_h^{\Omega\setminus S_{\delta}} (r)$ is a linear combination of $h\ZZ^d$ shifts of $\phi_{d,k}$
 taken from the set
 $$
 \left(h\ZZ^d\cap (\Omega\setminus S_{\delta}) \right)+ h \: \supp \mathbf{c}
 \subset
\Omega.
$$
Because $W\subset \Omega$ as well, the theorem follows.
\eop

It is perhaps worth noting that a minor modification of the previous proof reveals that if $\Omega\subset \RR^d$ is compact
and $g\in \con(\Omega)$ vanishes on $\partial \Omega$, then $g$ is in the uniform closure of $\spam_{z\in\Omega} \phi(\cdot - z)$. 

Furthermore, if we consider $f\in \con(\Omega)$ (not necessarily vanishing on $\partial \Omega$) and we choose $\delta, \epsilon>0$
then there exists a finite set $Y\subset \Omega_{\delta}$ (the $\delta$ neighborhood of $\Omega$) and coefficients $(a_y)_{y\in Y}\in \RR^Y$
so that $\|f - \sum_{y\in Y} a_y \phi_{d,k}(\cdot -y)\|_{\infty}<\epsilon$. This is done simply by considering a continuous extension of 
$f$ to $\RR^d$ that  vanishes on $\mathrm{cl}(\RR^d \setminus \Omega_{\delta})$ and applying the observation of the previous paragraph.

\section{Density in $\con(\ball^d)$}
In this section, we apply the scheme developed in Section 2 to prove fundamentality in $\Omega = \B^d$.
We will first show that we can approximate arbitrary continuous functions on the sphere using linear combinations of polyharmonic splines centered on the sphere. This fact, in conjunction with
Theorem~\ref{main_general} will deliver our density result.
\subsection{Approximation on the unit sphere}
We begin by introducing the Fourier-Gegenbauer expansion of a function $\psi:[-1,1]
\rightarrow \RR$. The Gegenbauer polynomials $P_j^{(\lambda)}$, $\lambda
> -1/2$, $j = 0,1,\ldots$, are orthonormal polynomials (we will assume this normalization 
in what follows) on the interval $(-1,1)$ with respect to the weight
$w^{(\lambda)}(t)=(1-t^2)^{\lambda-1/2}$. In other words,
\begin{equation}\label{orthonormality}
\int_{-1}^1 P_\ell^{(\lambda)} (t) P_j^{(\lambda)}(t) w^{(\lambda)}(t) \dif t = \delta_{\ell,j}, \quad \ell,j=0,1,\ldots .
\end{equation}
The Fourier-Gegenbauer coefficients of the function $\psi:[-1,1]\rightarrow \RR$ are
\beqn\label{geg_coeffs}
\alpha_j^\lambda(\psi) & = & \int_{-1}^1 \psi(t) P_j^{(\lambda)} (t)w^{(\lambda)}(t) \dif t.
\eeqn

The first step of our proof is to approximate $f \in \con(\ball^d)$
on the boundary $\sphere^{d-1}$ using polyharmonic splines. 
If $x,y\in \sphere^{d-1}$ then 
$|x|=|y|=1$ 
and we can use the definition (\ref{PHS_definition}) to rewrite  
$\phi_{d,k}(| x-y|)=\psi_{d,k}(x\cdot y)$, 
with
$$
\psi_{d,k}(r)= 
\begin{cases}
 (2-2r)^{k-d/2} , & d \;\; {\rm odd}, \\
\frac12 (2-2r)^{k-d/2} \log (2-2r), & d \;\; {\rm even}.
\end{cases}
$$
A crucial result for us is found in Sun and Cheney \cite[Theorem 2]{sun}.
\begin{proposition} \label{sun}
Fix $d$ and let $\lambda=(d-2)/2$. Let $\psi \in \con[-1,1]$. In order
that the set $\{ x \mapsto \psi(x\cdot y) : y \in \sphere^{d-1}\}$ be
fundamental in $\con(\sphere^{d-1})$ it is necessary and sufficient
that $\alpha_j^\lambda(\psi) \neq 0$, $j=0,1,\ldots$.
\end{proposition}

To simplify this, we use the families of functions $G_{\beta}(r) := \frac12 (2-2r)^{\beta} \log(2-2r)$
and $F_{\beta} (r) := (2-2r)^{\beta}$ 
defined for $\Re \beta>0$. We have
$\psi_{d,k}(r)= F_{k-d/2}(r)$ or $G_{k-d/2}(r)$, depending on the parity of $d$.

In \cite{BaHu}, the Fourier-Gegenbauer coefficients 
$a_j^\lambda(F_{\beta})$ and $a_j^{\lambda}(G_{\beta})$
are computed for  $\beta>0$ (in addition to coefficients for many other spherical basis functions $\psi$). 
We give a streamlined calculation below in (\ref{potential}) and (\ref{surface_spline}) respectively. 
For our immediate goal it is not important what the values
are, only that they are nonzero (see below for more details about these coefficients).

\begin{lemma}
For $d \ge 2$ and $k > d/2$, the Fourier-Gegenbauer coefficients of $\psi_{d,k}$ are all nonzero.
\end{lemma}

An immediate consequence of the above lemma, using Proposition \ref{sun}, we get the first new result presented in this paper:
\begin{theorem} \label{sunres}
The set $\{ \phi_{d,k}(| \cdot -y|) : y \in \sphere^{d-1}\}$ is dense in $\con (\sphere^{d-1})$.
\end{theorem}

This, in conjunction with Theorem \ref{main_general} proves the following.

\begin{theorem}\label{main}
Let $f \in \con(\ball^d)$. 
Given $\epsilon > 0$ 
there is a finite subset 
$Y \subset \ball^d$ 
and a set of coefficients 
$\{ \alpha_y : y \in Y \}$ 
such that
$$\left|f(x) - \sum_{y \in Y} \alpha_y \phi_{d,k}(|x-y|)\right| \le \epsilon$$ 
for all $x \in \ball^d$.
\end{theorem}

\subsection{Approximation on balls of arbitrary radius}\label{arbitrary_radius}
We now consider density of $\phi_{d,k}$ in $\ball^d(\rho)$, the ball of radius $\rho>0$. 
In this case, the proof of  Theorem~\ref{main_general} applies, if the analog
to Theorem~\ref{sunres} is established for the sphere of radius $\rho>0$: 
$\sphere^{d-1}(\rho) :=\{x\in \RR^d\mid |x| = \rho\}.$

Equivalently, this problem can be considered by simply rescaling $\phi_{d,k}  (\rho |x|)= \rho^{2k-d}\rbf(|x|)$ 
and considering the density of $\rbf$ on $\ball^d$ -- clearly, the translates of $\phi_{d,k}$ are dense 
in the set $\ball^{d}(\rho)$ if and only if
the translates of $\rbf$ are dense in $\ball^{d}$. A quick computation shows that
$$\rbf(r) = \begin{cases} r^{2k-d}, &\qquad d\text{ is odd}\\
r^{2k-d} \log r+ (\log \rho)\, r^{2k-d}, &\qquad d\text{ is even}.
\end{cases}$$
In odd dimensions, $\phi_{d,k} = \rbf$ is scale invariant, 
while in even dimensions it is necessary to 
add a polynomial term (in this case, $\rbf$ depends on $\rho$ as well as $d$ and $k$). 
Thus, the odd dimensional case follows directly from Theorem~\ref{main} and we are left to consider $d\in2\NN$.
The operator $\dlaplace^k$ annihilates polynomials of degree $2k-1$, so the only hinderance
to obtaining density in $\ball^d(\rho)$, is the fundamentality of
$\rbf(|x-y|)$ 
on $\sphere^{d-1}$.

Fix $d\in 2\NN$. The restriction of $\rbf$ to the unit sphere is $\rbf(|x-y|) = \widetilde{\psi}_{d,k}(x\cdot y)$. Indeed, 
\begin{eqnarray*}
\rbf(|x - y|)  &=& 
\frac12 (2-2x\cdot y)^{k-d/2} \log(2-2 x\cdot y)\\
&\mbox{}&+
 \log(\rho) (2-2x\cdot y)^{k-d/2}\\
&=&\psi_{d,k}(x\cdot y)+(\log\rho) p_{d,k}(x\cdot y),
\end{eqnarray*}
that is, $ \widetilde{\psi}_{d,k} = \psi_{d,k}+(\log\rho)\, p_{d,k}$.
Because $d$ is even, $ p_{d,k}(t)=(2-2t)^{k-d/2}$ 
is a polynomial of degree $k-d/2$. Thus, the coefficients of $\widetilde{\psi}_{d,k}$ for $j>k-d/2$ coincide with those of $\psi_{d,k}$ (which are nonzero), and we wish to calculate the coefficients $a_j^{\lambda}( \widetilde{\psi}_{d,k} )$ for $j\le k-d/2$.

We are concerned with  coefficients 
$a_j^{\lambda}(p_{d,k}) = a_j^{\lambda}(F_{k-d/2})$ and $a_j^{\lambda}(\psi_{d,k}) =a_j^{\lambda}(G_{k-d/2})$ for $j\le k-d/2$.
(Recall that $F_{\beta}$ 
and $G_{\beta}$ 
are functions in $\con([-1,1])$ whose Fourier-Gegenbauer coefficients have been determined for $\beta>0$ in \cite{BaHu}.)
 Specifically,  
 we wish to find radii $\rho$ so that $(\log \rho)a_j^{\lambda}(F_{k-d/2}) = - a_j^{\lambda}(G_{k-d/2})$.
 This would imply
 $$
a_j^{\lambda}(\psi_{d,k})+(\log\rho) a_j^{\lambda}(p_{d,k}) = 0 
$$
which, in turn, would indicate that $A(\ball^d(\rho),\phi)$ is not fundamental in $\con\bigl(\ball^d(\rho)\bigr)$, by Proposition \ref{sun}.

\paragraph{Fourier-Gegenbauer coefficients} Because we use Gegenbauer polynomials that are normalized in a slightly differently manner 
than those of \cite{BaHu},%
\footnote{Our Gegenbauer polynomials satisfy the orthonormality condition (\ref{orthonormality}), while those of Baxter and Hubbert have the standard, albeit more complicated, normalization 
$\int_{-1}^1|P_j^{(\lambda)}(t)|^2 w^{(\lambda)}(t) \dif t =d_j$ with $d_j$ given in \cite[(2.4)]{BaHu}, and  \cite[(22.2.3)]{AS}.}
we briefly repeat their argument here, with the necessary (minor) changes.
Let us first recall Rodrigues's formula,  \cite[(22.11.2)]{AS}, for Gegenbauer polynomials
$$
P_j^{(\lambda)}(t) = 
c_{j,\lambda} (1-t^2)^{1/2-\lambda} \frac{d^j}{dt^j}(1-t^2)^{j+\lambda-1/2}.
$$
The constant $c_{j,\lambda}$, when $j$ and $\lambda$  are not both $0$, is given by
$$
c_{j,\lambda} =
\frac{2^{j+\lambda}(-1)^j}{ \sqrt{2\pi}}\sqrt{j+\lambda}
\frac{\Gamma(j+\lambda)\sqrt{\Gamma(j+2\lambda)}}{\Gamma(2j+2\lambda)\sqrt{\Gamma(j+1)}}.
$$
When $\lambda$ and $j$ are both $0$ we have
$
c_{0,0} = 
1/{\sqrt{\pi}}.
$

In general, for a function $\psi\in C^{\infty}(-1,1)$ with suitable behavior at the boundary, 
we use (\ref{geg_coeffs}) in conjunction with Rodrigues's formula to obtain
\begin{eqnarray*}
a_j^{\lambda}(\psi) 
&=&
c_{j,\lambda} \int_{-1}^1 \psi(t) \frac{d^j}{dt^j}(1-t^2)^{j+\lambda-1/2} \dif t\\
&=&
(-1)^jc_{j,\lambda} \int_{-1}^1 \psi^{(j)}(t) (1-t^2)^{j+\lambda-1/2} \dif t.
\end{eqnarray*}
In the second line we have used integration by parts.

\paragraph{Coefficients for $F_{\beta}$.} 
We first compute coefficients for $F_{\beta}(t) =(2-2t)^{\beta}$, following \cite[(2.12)]{BaHu}. We note that 
$F_{\beta}^{(j)}(t)  = (-1)^j 2^{\beta} \frac{\Gamma(\beta+1)}{\Gamma(\beta-j+1)}(1-t)^{\beta-j} $.
This gives, for $j\le \beta$,
\begin{equation*}
a_j^{\lambda}
(F_{\beta})
=
c_{j,\lambda}  
2^{\beta}
\frac{\Gamma(\beta+1)}{\Gamma(\beta-j+1)}
\int_{-1}^1 (1-t)^{\beta+\lambda-1/2} (1+t)^{j+\lambda-1/2} \dif t.
\end{equation*}
Applying the change of variable $T = (t+1)/2$, we obtain
\begin{eqnarray*}
a_j^{\lambda}(F_{\beta})
&=&
c_{j,\lambda}  
2^{2\beta+j}
\frac{\Gamma(\beta+1)}{\Gamma(\beta-j+1)}
\int_{0}^1 (1-T)^{\beta+\lambda-1/2} T^{j+\lambda-1/2} \dif T\\
&=&
c_{j,\lambda}  
2^{2\beta+j}
\frac{\Gamma(\beta+1)}{\Gamma(\beta-j+1)}
\mathrm{B}(\beta+\lambda+1/2,j+\lambda+1/2)
\end{eqnarray*}
By applying the familiar identity for the beta function $\mathrm{B}(x,y) = \frac{\Gamma(x)\Gamma(y)}{\Gamma(x+y)}$,
and reorganizing the above expression, we obtain the formula obtained in \cite[(2.12)]{BaHu}:
\begin{equation}\label{potential}
a_j^{\lambda}(F_{\beta}) =
c_{j,\lambda}
2^j 
  \Gamma(j+\lambda+1/2)
\frac{2^{2\beta}\Gamma(\beta+1)\Gamma(\beta+\lambda+1/2)}{\Gamma(\beta-j+1)\Gamma(\beta +j+2\lambda+1)}.
\end{equation}
We note that this holds for all integers $j$, with the understanding that when 
$\beta$ is an integer and
$j\ge \beta +1$
then
$a_j^{\lambda}(F_{\beta}) =0$. It is not difficult to see that $\beta\mapsto a_j^{\lambda}(F_{\beta})$ is analytic in the right half plane $\Re\beta>0$.

\paragraph{Coefficients for $G_{\beta}$.}We now focus on $G_{\beta}(t) =\frac12 (2-2t)^{\beta}\log (2-2t)$. These coefficients are computed in  \cite[(2.19)]{BaHu}. We repeat this here for completeness.
We note that 
$$
\frac12 (2-2t)^{\beta}\log (2-2t) = \frac12 \frac{\partial}{\partial \beta} (2-2t)^{\beta},
$$
so we obtain the Fourier-Gegenbauer coefficients by taking derivatives with respect to $\beta$ in (\ref{potential}). In this case, we have
\begin{equation}\label{surface_spline}
a_j^{\lambda}(G_{\beta}) = \frac12 
c_{j,\lambda}
2^j 
  \Gamma(j+\lambda+1/2)
   \frac{\partial}{\partial \beta}
   \left(
   \frac{2^{2\beta}\Gamma(\beta+1)\Gamma(\beta+\lambda+1/2)}{\Gamma(\beta-j+1)\Gamma(\beta +j+2\lambda+1)}
   \right).
\end{equation}
Using the {\em digamma} function $\Psi(t) = \frac{\Gamma'(t)}{\Gamma(t)}$, we have
$$ a_j^{\lambda}(G_{\beta}) = 
\frac12
a_j^{\lambda}(F_{\beta}) \tau(\beta,j,\lambda)
$$
where
$$ 
\tau(\beta,j,\lambda)= 
\log 4 +\Psi(\beta+1) +\Psi(\beta+\lambda+1/2) - \Psi(\beta-j+1) -\Psi(\beta+j+2\lambda+1).
$$
For integers $j$, formulas \cite[(6.3.2)]{AS} and \cite[(6.3.4)]{AS}  indicate that
$$\Psi(j+1) = -\gamma +\sum_{n=1}^j \frac{1}{n} 
\quad \text{and} \quad
\Psi(j+1/2) = -\gamma-2\log 2 +\sum_{n=1}^j\frac{1}{n-1/2}.$$

It follows that when $\rho>0$ is chosen so that 
$\log \rho = - \frac12 \tau(k-d/2,j,\lambda)$ for some $j\le k-d/2$,
then 
$A(\sphere^{d-1}(\rho),\phi_{d,k})$ is not fundamental in $\sphere^{d-1}(\rho)$.

\begin{example}  We consider the thin plate spline example $k=2,d=2$. In this case, $k-d/2 =1$ and $\lambda=0$, 
so degenerate
examples occur when $\log \rho  = -\frac12 \tau(1,0,0) $ or $-\frac12 \tau(1,1,0)$. A simple computation
shows that 
$$\tau(1,0,0) = 1+2 -1-1 =1$$
and 
$$\tau(1,1,0) =  1+2 -0 - 3/2 =  \frac{3}{2} $$
from which we obtain $\rho = e^{-1/2}$ and $\rho = e^{-3/4}.$
\end{example}

\section{Conclusion}
The results of this paper are a first attempt at treating a challenging problem in RBF theory:
\begin{equation}\label{density_question}
\text{Find radial functions $\phi$  so that for all}\ \rho>0, \ 
A\bigl( \ball^d(\rho),\phi\bigr)
\text{ is fundamental in $\con\bigl(\ball^d(\rho)\bigr)$.}
\end{equation}
We note that when $\phi$ is homogeneous the problem reduces to determining if
$A\bigl( \ball^d,\phi\bigr)$ is fundamental in $\con\bigl(\ball^d\bigr)$.

This can be viewed as a modification of the density problem considered by Pinkus  
to characterize the radial functions $\phi$  so that  
$A(\RR^d,\phi)$ is fundamental in $\con(\ball^d(\rho))$ for all $\rho>0$.
In \cite{Pinkus} it was demonstrated that this holds for any radial $\phi\in L_1(\RR^d)\cap \con(\RR^d)$,
although the class of functions satisfying Pinkus's density problem is known to be significantly larger.
%

We consider the approach to solving (\ref{density_question}) suggested in this paper
in a more general light. Of the two components of the approximant -- the part that approximates
in $\con_0(\Omega)$ and the part that approximates in $\con(\partial \Omega)$ -- we have little
to say about the former, since not all 
radial basis functions are known to have a dilation invariant, rapidly decaying 
localization (this is a perk of working with polyharmonic splines). The latter component can
be generalized quite easily, however.

For a radial function $\Phi\in\con(\RR^d)$, the restriction of $(x,y)\mapsto \Phi(x-y)$ to $\sphere^{d-1}$
satisfies the relationship $\Phi(x-y)  = \psi(x\cdot y)$ for some univariate function $\psi:[-1,1]\to \RR$. 
The Fourier-Gegenbauer coefficients of $\psi$ are
tied to the distributional Fourier transform $\widehat{\Phi}$ of $\Phi$.  
If $\Phi\in L_1(\RR^d)$,
then $\widehat{\Phi} (\xi)= \sigma(|\xi|)$ for some $\sigma\in\con\bigl([0,\infty) \bigr)$.
If $\widehat{\Phi}$ is also in $L_1(\RR^d)$, 
then for each $j$, the coefficient
$a_j^{\lambda}(\psi) = \int_{-1}^1 \psi(t) P_j^{(\lambda)}(t) w^{(\lambda)}(t) \dif t$ 
is determined by the Weber-Schafheitlin type integral: 
\begin{equation}\label{W-S}
\alpha_j^{\lambda}(\psi) = \int_0^{\infty} t \sigma(t) \bigl(J_{j+\lambda}(t)\bigr)^2 \dif t
\end{equation} 
(see \cite[Proposition 1]{CasFil} for a proof, or the earlier version \cite{NW}). 
Here $\lambda = \frac{d-2}{2}$  
and $J_{\nu}$ is a  Bessel function of the first kind. It follows that if  
$\int_0^{\infty} t \sigma(t) \bigl(J_{j+\lambda}(t)\bigr)^2 \dif t\ne 0$ 
for all $j$ then 
$A(\sphere^{d-1},\Phi)$ 
is fundamental in
$\con(\sphere^{d-1})$. (This holds, for instance, if $\sigma(t)\ge0$ and is nonzero on a set of positive measure.) 

Recently, zu Castell and Filbir 
and, independently, Narcowich, Sun and Ward 
have extended this idea to radial, continuous, conditionally positive definite functions.
A function $\Phi\in \con(\RR^d)$ is conditionally positive definite (CPD)  of order $k$
if for each finite subset $\Xi\subset \RR^d$ the quadratic form 
$\mathbf{a}=(a_{\xi})_{\xi\in\Xi}\mapsto \sum_{\xi\in\Xi} \sum_{\zeta\in\Xi} \overline{a_{\zeta}}a_{\xi} \Phi(\xi -\zeta)$
is positive semi-definite for all vectors $\mathbf{a}$ satisfying $\sum_{\xi\in\Xi} a_{\xi} p(\xi)=0$ for all polynomials $p\in\Pi_{k-1}(\RR^d)$.
This class includes many radial basis functions that exhibit polynomial growth -- in particular the polyharmonic spline $\phi_{d,k}$ is conditionally positive definite of order $k$.

In \cite[Chapter II, Section 4.4]{GV} 
 Gelfand and Vilenkin characterize  
 such functions
 by way of the distributional Fourier transform $\widehat{\Phi}$ .
For any test function
$\gamma$ asymptotically  like $\mathcal{O}(|\xi|^{2k+1})$ at the origin,
the duality pairing
$\langle \widehat{\Phi},\gamma\rangle $ 
is given by 
$\langle \widehat{\Phi},\gamma\rangle 
=
\int_{\RR^d\setminus\{0\}}\gamma(\xi)\, \dif \mu(\xi)$, 
where $\mu$ is a positive Borel measure on $\RR^d\setminus\{0\}$
satisfying $\int_{\ball^d\setminus\{0\}} |\xi|^{2k} \dif \mu(\xi)<\infty$ and $\mu(\RR^d\setminus\ball^d) <\infty$.\footnote{
This can be improved by taking test functions $\gamma$ that vanish  like $\mathcal{O}(|\xi|^{2k})$ at the origin, in which case
we get the {\em generalized Fourier transform of order $k$}. In this case,
$\langle \widehat{\Phi},\gamma\rangle =\sum_{|\alpha|=2k}a_{\alpha} D^{\alpha} \gamma(0)+
\int_{\RR^d}\gamma(\xi)\, \dif \mu(\xi)$, where $(-1)^k \sum_{|\alpha|=2k}a_{\alpha} x^{\alpha} \in \Pi_{2k}(\RR^d)$
is a conditionally positive definite polynomial of order $k$.}

When $\Phi$ is radial, the distributional Fourier transform $\widehat{\Phi}$ 
is necessarily rotationally symmetric. 
It follows that there is a Borel measure
$\sigma$ on $(0,\infty)$  for which 
$\int_{\RR^d\setminus\{0\}} \gamma(\xi) \dif \mu(\xi) = \int_ 0^{\infty} \int_{\omega\in\sphere^{d-1}}
\gamma(t\omega)t^{d-1}\, \dif \omega \, \dif\sigma(t)$.
In \cite[(4.5)]{CasFil} and in \cite[Proposition 3.1]{NSW} a similar relationship to (\ref{W-S}) is shown for radial, continuous,
conditionally positive definite functions of order $k$. 
Namely, for $j> k$ the Fourier--Gegenbauer
coefficients are 
$$a_j^{\lambda}(\psi)  = \int_0^{\infty} t \left(J_{j+\lambda}(t)\right)^2 \dif \sigma(t).$$
We note that the zeros of $J_{j+\lambda}^2$ form a discrete set. It follows that 
if the (topological) support of $\sigma$ contains an open interval 
(e.g., if $\sigma$ is  absolutely continuous with respect to Lebesgue measure and 
$\dif \sigma(x) = f(x) \dif x$ and $f:\RR_+\to\RR_+$ is positive on an open interval) 
then for each $j>k$, 
the coefficient $a_j^{\lambda}(\psi)$ is positive.
By Proposition \ref{sun},  the set of rotations of $\psi$ modulo a polynomial $p_k\in \Pi_k(\RR)$, 
$\{x\mapsto \psi(x\cdot y) +p_k(x\cdot y)\mid y\in \sphere^{d-1}\}$, is
fundamental in $\sphere^{d-1}$. 
There is a unique radial polynomial $P_{2k}\in\Pi_{2k}(\RR^d)$ so that $p_k(x\cdot y) = P_{2k}(x-y)$ for $x,y\in\sphere^{d-1}$, and hence, modulo addition of a polynomial of degree $2k$, the restriction of $\Phi$ is fundamental in $\sphere^{d-1}$. 

This also permits the consideration of arbitrary spheres.
The problem of fundamentality of the radial, continuous, CPD function $\Phi$ on $\sphere^{d-1}(\rho)$ is equivalent to fundamentality of $\Phi_{\rho} := \Phi(\rho \cdot)$ on $\sphere^{d-1}$. Let $\psi_{\rho}\in C([-1,1])$ be the zonal function $\psi_{\rho} (x\cdot y) = \Phi_{\rho}(x-y)$.
Because the dilation by $\rho$ causes an inverse dilation
of its distributional Fourier transform, 
$\langle \Phi_{\rho},\widehat{\gamma}\rangle 
= 
\rho^{-d} \langle \Phi,\widehat{\gamma}(\cdot/\rho)\rangle
= 
\langle \widehat{\Phi},\gamma(\rho \cdot)\rangle,$
it follows that $\widehat{ \Phi_{\rho}}$ is represented by the measure $\sigma_{\rho}(E) = \rho^{1-d} \sigma(E/\rho)$.
Thus, for $j>k$ we have
$a_j^{\lambda}(\psi_\rho)  =  \int_0^{\infty} t \left(J_{j+\lambda}(t)\right)^2 \dif \sigma_{\rho}(t)
 = \rho^{2-d} \int_0^{\infty} t \left(J_{j+\lambda}(\rho t)\right)^2 \dif \sigma (t)
 .$
It again follows that, modulo the addition of a polynomial of degree $2k$ (depending on $\rho$),  $\phi$ is fundamental 
in $\sphere^{d-1}(\rho)$.

This is more or less what happens in Section 3.2, where the two conditionally positive
definite functions $\phi_{d,k}$ and $\rbf$ have the same generalized Fourier transform $\xi \mapsto |\xi|^{-2k}$,
but their Fourier--Gegenbauer coefficients disagree on lower frequencies.
\subsection*{Acknowledgment} The authors would like to thank Michael Johnson for a number of very interesting conversations about the subject matter of this paper. In addition, the referees provided a number of excellent comments and suggestions which, from the authors' points of view, improved the preceding exposition.
\bibliographystyle{siam}
\bibliography{HLevesley}

\begin{thebibliography}{10}

\bibitem{AS}
{\sc M.~Abramowitz and I.~A. Stegun}, {\em Handbook of mathematical functions
  with formulas, graphs, and mathematical tables}, vol.~55 of National Bureau
  of Standards Applied Mathematics Series, For sale by the Superintendent of
  Documents, U.S. Government Printing Office, Washington, D.C., 1964.

\bibitem{BSW}
{\sc K.~Ball, N.~Sivakumar, and J.~D. Ward}, {\em On the sensitivity of radial
  basis interpolation to minimal data separation distance}, Constr. Approx., 8
  (1992), pp.~401--426.

\bibitem{BaHu}
{\sc B.~J.~C. Baxter and S.~Hubbert}, {\em Radial basis functions for the
  sphere}, in Recent progress in multivariate approximation
  ({W}itten-{B}ommerholz, 2000), vol.~137 of Internat. Ser. Numer. Math.,
  Birkh\"auser, Basel, 2001, pp.~33--47.

\bibitem{brown}
{\sc A.~Brown}, {\em Uniform approximation by radial basis functions, {\em in
  advances in numerical analysis. {v}ol. {ii}}}, in Proceedings of the {F}ourth
  {S}ummer {S}chool in {N}umerical {A}nalysis held at the {U}niversity of
  {L}ancaster, {L}ancaster, {J}uly 15--{A}ugust 3, 1990, W.~Light, ed., Oxford
  Science Publications, The Clarendon Press Oxford University Press, New York,
  1992, pp.~viii+210.
\newblock Wavelets, subdivision algorithms, and radial basis functions.

\bibitem{cheney}
{\sc W.~Cheney and W.~Light}, {\em A course in approximation theory}, vol.~101
  of Graduate Studies in Mathematics, American Mathematical Society,
  Providence, RI, 2009.
\newblock Reprint of the 2000 original.

\bibitem{dBR}
{\sc C.~de~Boor and A.~Ron}, {\em Fourier analysis of the approximation power
  of principal shift-invariant spaces}, Constr. Approx., 8 (1992),
  pp.~427--462.

\bibitem{DynLevinRippa}
{\sc N.~Dyn, D.~Levin, and S.~Rippa}, {\em Numerical procedures for surface
  fitting of scattered data by radial functions}, SIAM J. Sci. Statist.
  Comput., 7 (1986), pp.~639--659.

\bibitem{GV}
{\sc I.~M. Gelfand and N.~Y. Vilenkin}, {\em Generalized functions. {V}ol. 4},
  Academic Press [Harcourt Brace Jovanovich Publishers], New York, 1964 [1977].
\newblock Applications of harmonic analysis, Translated from the Russian by
  Amiel Feinstein.

\bibitem{HangThesis}
{\sc T.~Hangelbroek}, {\em Approximation by scattered translates of the
  fundamental solution of the biharmonic equation on bounded domains}, 2007.
\newblock Thesis (Ph.D.)--The University of Wisconsin - Madison.

\bibitem{Hang}
\leavevmode\vrule height 2pt depth -1.6pt width 23pt, {\em Error estimates for
  thin plate spline approximation in the disk}, Constr. Approx., 28 (2008),
  pp.~27--59.

\bibitem{Jackson}
{\sc I.~R.~H. Jackson}, {\em Convergence properties of radial basis functions},
  Constr. Approx., 4 (1988), pp.~243--264.

\bibitem{JiaMic}
{\sc R.~Q. Jia and C.~A. Micchelli}, {\em Using the refinement equations for
  the construction of pre-wavelets. {II}. {P}owers of two}, in Curves and
  surfaces ({C}hamonix-{M}ont-{B}lanc, 1990), Academic Press, Boston, MA, 1991,
  pp.~209--246.

\bibitem{johnson}
{\sc M.~Johnson}, {\em Private communication}.
\newblock to J. Levesley, 1999.

\bibitem{light2}
{\sc W.~Light and H.~Wayne}, {\em Spaces of distributions, interpolation by
  translates of a basis function and error estimates}, Numer. Math., 81 (1999),
  pp.~415--450.

\bibitem{NSW}
{\sc F.~J. Narcowich, X.~Sun, and J.~D. Ward}, {\em Approximation power of
  {RBF}s and their associated {SBF}s: a connection}, Adv. Comput. Math., 27
  (2007), pp.~107--124.

\bibitem{NW}
{\sc F.~J. Narcowich and J.~D. Ward}, {\em Scattered data interpolation on
  spheres: error estimates and locally supported basis functions}, SIAM J.
  Math. Anal., 33 (2002), pp.~1393--1410.

\bibitem{Pinkus}
{\sc A.~Pinkus}, {\em Some density problems in multivariate approximation}, in
  Approximation theory ({W}itten, 1995), vol.~86 of Math. Res.,
  Akademie-Verlag, Berlin, 1995, pp.~277--284.

\bibitem{CR1}
{\sc C.~Rabut}, {\em Elementary {$m$}-harmonic cardinal {$B$}-splines}, Numer.
  Algorithms, 2 (1992), pp.~39--61.

\bibitem{sun}
{\sc X.~Sun and E.~W. Cheney}, {\em Fundamental sets of continuous functions on
  spheres}, Constr. Approx., 13 (1997), pp.~245--250.

\bibitem{CasFil}
{\sc W.~zu~Castell and F.~Filbir}, {\em Radial basis functions and
  corresponding zonal series expansions on the sphere}, J. Approx. Theory, 134
  (2005), pp.~65--79.

\end{thebibliography}

\end{document}